\documentclass[12pt]{amsart}
\usepackage{amsmath, amssymb, euscript, hyperref}

\newtheorem{thm}{Theorem}[section]
\newtheorem{lem}[thm]{Lemma}
\newtheorem{prop}[thm]{Proposition}
\theoremstyle{definition}
\newtheorem{defn}[thm]{Definition}

\newtheorem{rem}[thm]{Remark}

\newcommand{\blackboard}[1]{\ensuremath{\mathbb{#1}}}

\newcommand{\Z}{\blackboard{Z}}

\begin{document}

\bibliographystyle{amsplain}

\address{Azer Akhmedov, Department of Mathematics,
North Dakota State University,
Fargo, ND, 58102, USA}
\email{azer.akhmedov@ndsu.edu}

\begin{center} {\bf On groups of diffeomorphisms of the interval with finitely many fixed points} \end{center}

\medskip

 \begin{center} Azer Akhmedov \end{center}

\bigskip

{\bf Abstract:} {\Small In \cite{N1}, it is proved that any subgroup of $\mathrm{Diff}_{+}^{\omega }(I)$ (the group of orientation preserving analytic diffeomorphisms of the interval) is either metabelian or does not satisfy a law. A stronger question is asked whether or not the Girth Alternative holds for subgroups of $\mathrm{Diff}_{+}^{\omega }(I)$. In this paper, we answer this question affirmatively for even a larger class of groups of orientation preserving diffeomorphisms of the interval where every non-identity element has finitely many fixed points. We show that every such (irreducible) group is either affine (in particular, metabelian) or has infinite girth. The proof is based on our study of discrete subgroups of the diffeomorphism group  $\mathrm{Diff}_{+}(I)$ which we initiated in \cite{A9} and later developed in  \cite{A1} and \cite{A2}; more specifically, our results are obtained by sharpening the tools from the earlier works \cite{A1} and \cite{A2}. One of the major tools (local transitivity)  is heavily exploited in \cite{A2} to get an extension of H\"older's theorem which is crucially used in this paper. We show that local transitivity can be proved for any (up to a conjugacy) non-affine group of  irreducible diffeomorphisms with every non-identity element having finitely many fixed points.} 

\vspace{1cm}

  \section{Introduction}
  
  Throughout this paper we will write $\mathbf{\Phi }$ (resp. $\mathbf{\Phi }^{\mathrm{diff}}$) to denote the class of subgroups of $\Gamma \leq \mathrm{Homeo}_{+}(I)$ (resp. $\Gamma \leq \mathrm{Diff}_{+}(I)$) such that every non-identity element of $\Gamma $ has finitely many fixed points. An important class of such groups is provided by $\mathrm{Diff}_{+}^{\omega }(I)$ - the group of orientation preserving analytic diffeomorphisms of $I$. Interestingly, not every group in $\mathbf{\Phi }$ is conjugate (or even isomorphic) to a subgroup of $\mathrm{Diff}_{+}^{\omega }(I)$. We will consider a natural metric on $\mathrm{Homeo}_{+}(I)$ induced by the $C^0$-metric by letting $||f|| = \displaystyle \mathop{\sup }_{x\in [0,1]}|f(x)-x|$.   
  
 \medskip
  
  For an integer $N\geq 0$ we will write $\mathbf{\Phi }_N$ (resp. $\mathbf{\Phi }^{\mathrm{diff}}_N$) to denote the class of subgroups of $\mathrm{Homeo}_{+}(I)$ (resp. $\mathrm{Diff}_{+}(I)$) where every non-identity element has at most $N$ fixed points. It has been proved in \cite{A1} that, for $N\geq 2$, any subgroup of $\mathbf{\Phi }^{\mathrm{diff}}_N$ of regularity $C^{1+\epsilon }$ is indeed solvable, moreover, in the regularity $C^2$ we can claim that it is metabelian. In \cite{A3}, we improve these results and give a complete classification of subgroups of $\mathbf{\Phi }^{\mathrm{diff}}_N, N\geq 2$, even at $C^1$-regularity. There, it is also pointed out that, in a strict sense, the presented classification picture fails in the continuous category, i.e. within the larger class $\mathbf{\Phi }_N$.

  \medskip
  
  One of the main results of this paper is the following 
  
  \medskip

 \begin{thm}\label{thm:main} Let $\Gamma \leq \mathrm{Diff}_{+}(I)$ be a finitely generated subgroup such that every non-identity element has finitely many fixed points. Then either $\Gamma $ is isomorphic to a subgroup of $\mathrm{Aff}_{+}(\mathbb{R})^n$ for some $n\geq 1$, or it has infinite girth.  
 \end{thm} 
  
 \medskip
 
  \begin{rem} Let us recall that the girth of a graph is the length of a shortest cycle in it (if there is no cycle, then the girth is equal to infinity), and the girth of a finitely generated group $\Gamma $ is defined as $$girth(\Gamma ) =\sup \{girth (\mathcal{G}(\Gamma , S)) : \langle S \rangle = \Gamma , |S| < \infty   \}$$ where $\mathcal{G}(\Gamma , S)$ denotes the Cayley graph of $\Gamma $ with respect to a finite generating set $S$. We refer the reader to \cite{S} or \cite{A6} for basic properties of girth, but would like to recall here the most relevant fact that  if a finitely generated group satisfies a law and is not isomorphic to $\Z $, then it has a finite girth. In particular, the girth of a finitely generated non-cyclic virtually solvable group is finite. Thus, as a corollary of Theorem  \ref{thm:main}, we obtain that a finitely generated subgroup of $\mathrm{Diff}_{+}(I)$ where every non-identity diffeomorphism has finitely many fixed points, satisfies no law unless it is isomorphic to a subgroup of $\mathrm{Aff}_{+}(\mathbb{R})^n$ for some $n\geq 1$. This also implies a positive answer to Question (v) from \cite{N1}. If $\Gamma $ is {\em irreducible} (i.e. it does not have a global fixed point in $(0,1)$), then one can take $n=1$.
  \end{rem}

  The proof of Theorem \ref{thm:main} relies on the study of diffeomorphism groups which act locally transitively. To be precise, we need the following definitions.
  
  \begin{defn} A subgroup $\Gamma \leq \mathrm{Diff}_{+}(I)$ is called {\em locally transitive} if for all $p\in (0,1)$ and $\epsilon > 0$, there exists $\gamma \in \Gamma $ such that $||\gamma || < \epsilon $ and $\gamma (p)\neq p$.  
  \end{defn}
  
  \medskip
  
  \begin{defn} A subgroup $\Gamma \leq \mathrm{Diff}_{+}(I)$ is called {\em dynamically 1-transitive} if for all $p\in (0,1)$ and for every non-empty open interval $J\subset (0,1)$ there exists $\gamma \in \Gamma $ such that $\gamma (p)\in J$.
  \end{defn}
  
  \medskip
  
  The notion of dynamical $k$-transitivity is introduced in \cite{A2} where we also make a simple observation that local transitivity implies dynamical 1-transitivity. Notice that if the group is dynamically $k$-transitive for an arbitrary $k\geq 1$,then it is dense in the $C^0$-metric. As we emphasize in\cite{A2} and \cite{A8}, dynamical $k$-transitivity for some high values of $k$ (even for $k=2$) is usually very hard, if possible, to achieve, and would be immensely useful (all the results obtained in \cite{A2} are based on just establishing the local transitivity which implies dynamical $k$-transitivity only for $k=1$). In this paper, we would like to introduce a notion of {\em weak transitivity} which turns out to be sufficient for our purposes but it is also interesting independently.
  
  \begin{defn} A subgroup $\Gamma \leq \mathrm{Diff}_{+}(I)$ is called {\em weakly $k$-transitive} if for all $g\in \Gamma $, for all $k$ points $p_1, \dots , p_k\in (0,1)$ with $p_1 < \dots < p_k$, and for all $\epsilon > 0$, there exist $\gamma \in \Gamma $ such that $g(\gamma (p_i))\in (\gamma (p_{i-1}), \gamma (p_{i+1}))$, for all $i\in \{1, \dots , k\}$ and $\gamma (p_k) < \epsilon $ where we assume $p_0 = 0, p_{k+1} = 1$. We say $\Gamma $ is {\em weakly transitive} if it is weakly $k$-transitive for all $k\geq 1$. 
  \end{defn}
  
  The proof of the Theorem \ref{thm:main} will follow immediately from the following four propositions which seem interesting to us independently.
  
  \begin{prop}\label{thm:first} Any irreducible subgroup of $\mathbf{\Phi }^{\mathrm{diff}}$ is either isomorphic to a subgroup of $\mathrm{Aff}_{+}(\mathbb{R})$ or it is locally transitive.
  \end{prop}
  
  \begin{prop} \label{thm:weakt} Any irreducible subgroup of $\mathbf{\Phi }^{\mathrm{diff}}$ is either isomorphic to a subgroup of $\mathrm{Aff}_{+}(\mathbb{R})$ or it is weakly transitive.
 \end{prop}
  
  \begin{prop}\label{thm:second} Any finitely generated locally transitive irreducible subgroup of $\mathbf{\Phi }^{\mathrm{diff}}$ is either isomorphic to a subgroup of $\mathrm{Aff}_{+}(\mathbb{R})$ or it has infinite girth.
  \end{prop}   
  
  \begin{prop}\label{thm:third} For any $N\geq 1$, any irreducible subgroup of $\mathrm{Diff}_{+}(I)$ where every non-identity element has at most $N$ fixed points is isomorphic to a subgroup of $\mathrm{Aff}_{+}(\mathbb{R})$.   
  \end{prop}
  
  \begin{rem}  \label{thm:remark}  Let us point out that, in $C^{1+\epsilon }$-regularity, Proposition \ref{thm:third} (the extension of H\"older's Theorem) is proved in \cite{A2} under somewhat weaker conclusion, namely, by replacing the condition ``isomorphic to a subgroup of $\mathrm{Aff}_{+}(\mathbb{R})$" with ``solvable", and in $C^2$-regularity by replacing it with ``metabelian". Subsequently, another (simpler) argument for this result is given in the analytic category, by A.Navas \cite{N1}. In the $C^1$ class, Proposition \ref{thm:third} is proved in \cite{A3} where it is shown that non-affine irreducible subgroups are non-discrete in the $C^0$ metric; this is sufficient for the claim of Proposition \ref{thm:third}.  {\em Alternatively}, let us recall that in \cite{A1}, we have indeed proved the following claim: {\em Let $\Gamma \leq \mathrm{Diff}_{+}(I)$ be an irreducible subgroup containing a free semigroup in two generators $f, g$ such that either $f'(0) = g'(0) = 1$ or $f'(1) = g'(1) = 1$. Then $\Gamma $ is not $C_0$-discrete, moreover, there exists non-identity elements in $[\Gamma , \Gamma ]$ arbitrarily close to the identity in the $C_0$ metric}. On the other hand, by Mueller-Tsuboi trick \cite{MN, T, M} we can always arrange the derivative at the end point 1 (or at the end point 0) to be 1, hence we obtain Proposition \ref{thm:third}, as well as the above claim about the  existence  non-identity elements in $[\Gamma , \Gamma ]$ arbitrarily close to the identity, independently of the discussion in \cite{A3}. 
  \end{rem}
  
  \medskip 
  
  \begin{rem} The major key tool of obtaining an extension of H\"older's Theorem in \cite{A2} is {\em local transitivity}. It implies  {\em dynamical 1-transitivity} (proximality) but it is significantly stronger than that. While the  dynamical 1-transitivity can be obtained in a more general setting of   $\mathbf{\Phi }$ (see the discussion in \cite{A3}, pp.11-13), local transitivity seems to be a more sensitive property; \cite{A3} makes no direct efforts in this direction even for some restricted sub-classes. Proposition \ref{thm:weakt} provides local transitivity property in the setting of the whole class  $\mathbf{\Phi }^{\mathrm{diff}}$ (instead of just $\mathbf{\Phi }^{\mathrm{diff}}_N$ as in \cite{A2}, see Proposition 1.2 and 1.4 there). We also need to emphasize that in the proof of Theorem \ref{thm:main} , we use only the dynamical 1-transitivity instead of full strength of local transitivity. 
  \end{rem} 
  
   Since Proposition \ref{thm:third} is established, we need to prove only Proposition \ref{thm:first}, \ref{thm:weakt} and \ref{thm:second}. The first of these propositions is proved in Section 2, by modifying the main argument from \cite{A1} (Let us also point out that another modification of the main result of \cite{A1} has been presented in \cite{A3}). The second and third propositions are proved in Section 3; in the proof of the third proposition, we follow the main idea from \cite{A4}, but our argument here is simpler.  Notice that Proposition   \ref{thm:second} establishes the so-called Girth Alternative in the class   $\mathbf{\Phi }^{\mathrm{diff}}$. In any class $\mathcal{C}$ of finitely generated groups, $\mathcal{C}$ is said to satisfy the \textit{Girth Alternative} if any group from the class $\mathcal{C}$ has either infinite girth or is virtually solvable. We refer the reader to \cite{S} or \cite{A6} for basic properties of girth; in \cite{A7}, the Girth Alternative is proved for the classes of word hyperbolic, one-relator and linear groups, and in \cite{A4}, it is proved for the class $PL_{+}(I)$.
   
  \medskip

   For all $f\in \mathrm{Homeo}_{+}(I)$ we will write $Fix(f) = \{x\in (0,1) \ | \ f(x) = x\}$. The group $\mathrm{Aff}_{+}(\mathbb{R})$ will denote the group of all orientation preserving affine homeomorphisms of $\mathbb{R}$, i.e. the maps of the form $f(x) = ax+b$ where $a>0$. The conjugation of this group by an orientation preserving homeomorphism $\phi : I \to \mathbb{R}$ to the group of homeomorphisms of the interval $I$ will be denoted by $\mathrm{Aff}_{+}(I)$ (we will drop the conjugating map $\phi $ from the notation). By choosing $\phi $ appropriately, one can also conjugate $\mathrm{Aff}_{+}(\mathbb{R})$ to a group of diffeomorhisms of the interval as well, and (by fixing $\phi $) we will refer to it as the group of affine diffeomorphisms of the interval.      
  
  \medskip
  
   Having the notion of conjugacy, we can indeed state a stronger claim, replacing ``isomorhism'' with actual ``conjugacy'' in our propositions. In the following theorem we make this stronger statement and also unite the claims of Propositions \ref{thm:first},    \ref{thm:weakt},  and \ref{thm:second}. 
   
   \begin{thm}\label{thm:main2} Any irreducible non-Abelian group  $\Gamma \in \mathbf{\Phi }^{\mathrm{diff}}$ is either conjugate to a subgroup of $\mathrm{Aff}_{+}(\mathbb{R})$ or it is locally transitive, weakly transitive, and if finitely generated, has infinite girth.    
   \end{thm}
   
   We are also referring the reader to Theorem 1.1 of \cite{A3}, where a strengthening of our fourth proposition (i.e. Proposition  \ref{thm:third}) is presented replacing again the isomorphism claim with a stronger claim of conjugacy. Let us remind that a conjugacy claim is also proved by Farb-Franks \cite{FF}  for non-Abelian subgroups $\Gamma \leq  \mathrm{Diff}_{+}^2(\mathbb{R}) $ in the case when every non-identity diffeomorphism has at most one fixed point (see Theorem 1.5 there). In the proof, the authors establish that there is no wandering interval of the form $\theta ^{-1}([x_0, x_1], x_0 < x_1$ where $\theta  = \nu [0, x)$ is the semi-conjugacy defined there. Our conjugacy claim in Theorem \ref{thm:main2} also follows from the non-existence of wandering intervals (as in \cite{FF}), but in our case, this claim (i.e. non-existence of wandering intervals) follows directly from the local transitivity.

\medskip 

 \begin{rem} Let us also emphasize that, in $\mathrm{Homeo}_{+}(\mathbb{R}) $, faithful representations of metabelian groups  such as $BS(1, m), m\geq 2$, are well-known  where every  element has at most $N$ fixed points for some uniform $N\geq 2$ and at least one element having exactly $N$ fixed points (these representations are even semi-conjugate to affine one). Theorem \ref{thm:main2} implies that none of these representations can exist in $\mathrm{Diff}_{+}(I)$.
 \end{rem}

  {\bf Notations:} We would like to end the introductory section by introducing some useful notations. If $\Gamma $ is a subgroup from the class $\mathbf{\Phi }$ then one can introduce the following biorder in $\Gamma $: for $f, g\in \Gamma $, we let $g < f$ if $g(x) < f(x)$ near zero. If $f$ is a positive element then we will also write $g<<f$ if $g^n<f$ for every integer $n$; we will say that $g$ is {\em infinitesimal} with respect to $f$. For $f\in \Gamma $, we write $\Gamma _f = \{\gamma \in \Gamma : \gamma << f \}$ (so $\Gamma _{f}$ consists of diffeomorphisms which are infinitesimal with respect to $f$). Notice that if $\Gamma $ is finitely generated with a fixed finite symmetric generating set, and $f$ is the biggest generator, then $\Gamma _f$ is a normal subgroup of $\Gamma $, moreover, $\Gamma /\Gamma _f$ is Archimedean, hence Abelian, and thus we also see that $[\Gamma , \Gamma ]\leq \Gamma _{f}$.

  \vspace{1cm}
  
  {\em Acknowledgement:} We are thankful to Michele Triestino for bringing the references \cite{MN, M} to our attention. 
  
   \vspace{1cm}

   \section{$C_0$-discrete subgroups of $\mathrm{Diff}_{+}(I)$: strengthening the results of \cite{A1}} 
   
   \medskip

 Let us first quote the following theorem from \cite{A1}. 
 
 \begin{thm}[Theorem A] Let $\Gamma \leq \mathrm{Diff}_{+}(I)$ be a subgroup such that $[\Gamma ,\Gamma ]$ contains a free semigroup on two generators. Then $\Gamma $ is not $C_0$-discrete, moreover, there exists non-identity elements in $[\Gamma , \Gamma ]$ arbitrarily close to the identity in $C_0$ metric.   
  \end{thm}

 In \cite{A2}, Theorem A is used to obtain local transitivity results in $C^{2}$-regularity ($C^{1+\epsilon }$-regularity ) for subgroups from $\mathbf{\Phi }^{\mathrm{diff}}_N$ which have derived length at least three (at least $k(\epsilon )$). However, we need to obtain local transitivity results for subgroups which are 1) non-abelian (not necessarily non-metabelian);  2) from a larger class $\mathbf{\Phi }^{\mathrm{diff}}$; 3) and have only $C^1$-regularity.
 
 \medskip
 
 To do this, first we observe that within the class $\mathbf{\Phi }^{\mathrm{diff}}$, the condition ``$\Gamma $ contains a free semigroup" by itself implies that ``$[\Gamma , \Gamma ]$ is either Abelian or contains a free semigroup".
         
 \medskip
 
 \begin{prop}\label{prop:commutator} Let $\Gamma $ be an irreducible non-Abelian subgroup from the class $\mathbf{\Phi }^{\mathrm{diff}}$. Then either $\Gamma $ is locally transitive or $[\Gamma , \Gamma ]$ contains a free semigroup on two generators. 
 \end{prop}
 
 \medskip
 
  For the proof of the proposition, we need the following 
  
 \begin{defn} Let $f, g\in \mathrm{Homeo}_{+}(I)$. We say the pair $(f,g)$ is {\em crossed} if there exists a non-empty open interval $(a, b)\subset (0,1)$ such that one of the homeomorphisms fixes $a$ and $b$ but no other point in $(a,b)$ while the other homeomorphism maps either $a$ or $b$ into $(a,b)$.
  \end{defn}
  
 \medskip
 
  It is a well known folklore result that if $(f,g)$ is a crossed pair then the subgroup generated by $f$ and $g$ contains a free semigroup on two generators (see \cite{N2}).

 {\bf Proof of Proposition \ref{prop:commutator}.}  Let us first assume that $\Gamma $ is metabelian, and let $N$ be a nontrivial Abelian normal subgroup of $\Gamma $. By H\"older's Theorem, there exists $f\in \Gamma $ such that $Fix(f)\neq \emptyset $. On the other hand, by irreducibility of $\Gamma $, $Fix(g) = \emptyset $ for all $g\in N\backslash \{1\}$.  
 
 \medskip
 
 Let us recall that (see Proposition 3.2. in \cite{A3}) any group from the class  $\mathbf{\Phi }$ is either Abelian or contains a free semigroup on two generators. Then by the Remark \ref{thm:remark} in the introduction, the subgroup $\langle f, N\rangle $ generated by $f$ and $N$ contains elements in its commutator subgroup (hence in $N$) arbitrarily close to the identity. Thus, for all $\epsilon > 0$, there exists $\omega \in N$ such that $||\omega || < \epsilon $. Since $N$ acts freely, this implies that $N$, hence $\Gamma $,  is locally transitive.\footnote{In the analytical category, non-discreteness of $N$ is also proved in Lemma 6.2 in \cite{A5}.}
 
 \medskip
  
 Let us now assume that $\Gamma $ has derived length more than 2 (possibly infinity). Then $\Gamma ^{(1)} = [\Gamma , \Gamma ]$ is not Abelian. Then by H\"older's Theorem there exists two elements $f, g\in \Gamma ^{(1)}$ such that $Fix(f) \neq Fix(g)$. We may assume that (by switching $f$ and $g$ if necessary) there exists $p, q\in Fix(f)\cup \{0,1\}$ such that $Fix(f)\cap (p,q) = \emptyset, Fix(g)\cap (p,q)\neq \emptyset $. Without loss of generality we may also assume that $f(x) > x, \forall x\in (p,q)$ and $g(p) \geq p$. If $g(p) > p$ then $f$ and $g$ form a crossed pair. But if $g(p) = p$ then for sufficiently big $n$, $f$ and $g^n$ form a crossed pair. $\square $

 \medskip
 
 To finish the proof of Proposition \ref{thm:first} it remains to show the following

 \begin{thm}\label{thm:A'} Let $\Gamma \leq \mathbf{\Phi }^{\mathrm{diff}}$ be an irreducible subgroup containing a free semigroup in two generators such that $f'(0) = f'(1) = 1$ for all $f\in \Gamma $. Then $\Gamma $ is locally transitive.  
  \end{thm}

  \medskip
  
 {\bf Proof.} The proof is very similar to the proof of Theorem A from \cite{A1} with a crucial extra detail (and with some simplifications).

  \medskip
  
  Without loss of generality, we may assume that $\Gamma $ is irreducible. Let $p\in (0,1), \epsilon > 0, M = \displaystyle \mathop{\sup }_{0\leq x\leq 1}(|f'(x)|+|g'(x)|)$. Let also $N\in \mathbb{N}$ and $\delta > 0$ such that $1/N < \min \{ \epsilon , p, 1-p\}$ and for all $x\in [1-\delta , 1)$, the inequalities $|\phi '(x)-1| < 1/10$ and $\phi (x) \neq x$ hold where $\phi \in \{f, g, f^{-1}, g^{-1}\}$.

 \medskip
 
 Let $W = W(f,g)$ be an element of $\Gamma $ such that $$(\{f^iW(1/N) \ | \ -2\leq i\leq 2\}\cup \{g^iW(1/N) \ | \ -2\leq i\leq 2\})\subset [1-\delta , 1]$$ and let $m$ be the length of the reduced word $W$. Let also $x_i = i/N, 0\leq i\leq N$ and $z = W(1/N)$.     
  
  \medskip
  
 By replacing one or both of $f, g$ with $f^{-1},g^{-1}$ respectively if necessary, we may assume that $f(z)> z$ and $g(z) > z$. Let $z_0 = z, \alpha = f $ and $\beta = g$ \ \ \footnote{we make these change of letters if a reader wants to compare our argument with the argument of \cite{A1}.} 
   
   \medskip
   
   Thus $\displaystyle \mathop{\sup }_{0\leq x\leq 1}(|\alpha '(x)|+|\beta '(x)|) =  M$ and the length of the word $W$ in the alphabet $\{\alpha , \beta , \alpha ^{-1}, \beta ^{-1}\}$ is $m$. 
  
  \medskip
 
 Now, for every $n\in \mathbb{N}$, let $$\mathbb{S}_n = \{U(\alpha ,\beta )  \ | \ U(\alpha ,\beta ) \ \mathrm{is \ a \ positive \ word \ in} \ \alpha , \beta  \ \mathrm{of \ length } \ n.\}$$ 
 
  \medskip

 Then $|\mathbb{S}_n| = 2^n$ and $V(z_0) \geq z_0$ for all $V\in \mathbb{S}_n$.
 
 \medskip
   
  Now, let $c_k = \frac{1}{100} - \frac{1}{100(k+1)}$ for all $k\geq 1$ (notice that the sequence $(2-c_k)_{k\geq 1}$ is decreasing to 1.99).  Then there exits a natural $n_0$ such that for all $n\geq n_0$ the following two conditions are satisfied:
 
 \medskip
 
 (i) there exist a subset $\mathbb{S}_n^{(1)} \subseteq \mathbb{S}_n$ such that $|\mathbb{S}_n^{(1)}| > (2-c_1)^n$, and for all $g_1, g_2\in \mathbb{S}_n^{(1)}$, $$|g_1W(x)-g_2W(x)| < \frac{1}{(1.9)^n}, \forall x\in \{x_i \ | \ 1\leq i\leq N-1\}\cup \{p\}.$$
 
 \medskip
 
 (ii) $M^{m+1}(1.1)^n\frac{1}{(1.9)^n} < \epsilon $.
 
 \bigskip
 
 Let us also fix  any $n$ greater than $n_0$. Then from Mean Value Theorem we obtain that $$|(g_1W)^{-1}(g_2W)(x) - x| < 2\epsilon $$ for all $x\in [0,1]$. \footnote{to explain this, we borrow the following computation from \cite{A1}: let $$h_1 = g_1W, h_2 = g_2W, y_i = W(x_i), z_i' = g_1(y_i), z_i'' = g_2(y_i), 1\leq i\leq N.$$ 
  
 \medskip
 
 Then for all $i\in \{1, \ldots , N-1\}$, we have 
 
 $$|h_1^{-1}h_2(x_i) - x_i| = |(g_1W)^{-1}(g_2W)(x_i) - x_i| =$$ \ $$|(g_1W)^{-1}(g_2W)(x_i) - (g_1W)^{-1}(g_1W)(x_i)| = |W^{-1}g_1^{-1}g_2(y_i) - W^{-1}g_1^{-1}g_1(y_i)|$$ \ $$ = |
 W^{-1}g_1^{-1}(z_i'') - W^{-1}g_1^{-1}(z_i')|$$
 
 \medskip
 
  Since $g_1, g_2\in \mathbb{S}_n$, by the Mean Value Theorem, we have $$|h_1^{-1}h_2(x_i) - x_i|\leq M^{m+1}(1.1)^n|z_1'-z_1''| < M^{m+1}(1.1)^n\frac{1}{(1.9)^n}$$ where the latter inequality $|z_1'-z_1''| < \frac{1}{(1.9)^n}$ is holding for all $g_1, g_2\in \mathbb{S}_n^{(1)}$  by our choice of the set $\mathbb{S}_n^{(1)}$. Since $m$ is fixed, for sufficiently big $n$ we obtain that $|h_1^{-1}h_2(x_i) - x_i| < \epsilon $. Then $|h_1^{-1}h_2(x) - x| < 2\epsilon $ for all $x\in [0,1]$.} 
  
  If $g_1W(p)\neq g_2W(p)$ for some $g_1, g_2\in \mathbb{S}_n^{(1)}$ then we are done. Otherwise we define the next set $$\mathbb{S}_n^{(2)} \subseteq \{U(\alpha ,\beta )g \ | \ g\in \mathbb{S}_n^{(1)}, U(\alpha ,\beta ) \ \mathrm{is \ a \ positive \ word \ in} \ \alpha , \beta  \ \mathrm{of \ length} \ n.\}$$  such that $|\mathbb{S}_n^{(2)}| > (2-c_2)^{2n}$, and for all $g_1, g_2\in \mathbb{S}_n^{(2)}$, $$|g_1W(x)-g_2W(x)| < \frac{1}{(1.9)^{2n}}, \forall x\in \{x_i \ | \ 1\leq i\leq N-1\}\cup \{p\}.$$

 \medskip 
 
 Again, by Mean Value Theorem, we obtain that $$|(g_1W)^{-1}(g_2W)(x) - x| < 2\epsilon $$ for all $g_1, g_2\in \mathbb{S}_n^{(2)}, x\in [0,1]$. Again, if $g_1W(p)\neq g_2W(p)$ for some $g_1, g_2\in \mathbb{S}_n^{(2)}$ then we are done, otherwise we continue the process as follows: if the sets $\mathbb{S}_n^{(1)}, \dots ,  \mathbb{S}_n^{(k)}$ are chosen and $g_1W(p) = g_2W(p)$ for all $g_1, g_2\in \mathbb{S}_n^{(k)}$, then we choose $$\mathbb{S}_n^{(k+1)} \subseteq \{U(\alpha ,\beta )g \ | \ g\in \mathbb{S}_n^{(k)}, U(\alpha ,\beta ) \ \mathrm{is \ a \ positive \ word \ in} \ \alpha , \beta  \ \mathrm{of \ length} \ n.\}$$  such that $|\mathbb{S}_n^{(k+1)}| > (2-c_{k+1})^{(k+1)n}$, and for all $g_1, g_2\in \mathbb{S}_n^{(k+1)}$, $$|g_1W(x)-g_2W(x)| < \frac{1}{(1.9)^{(k+1)n}}, \forall x\in \{x_i \ | \ 1\leq i\leq N-1\}\cup \{p\}.$$
 
 \medskip
 
 Since the sequence $(2-c_k)_{k\geq 1}$ is decreasing to 1.99, we can choose the same sufficiently big $n$ at each step (i.e. the same fixed $n$). But $\Gamma $ belongs to the class $\mathbf{\Phi }^{\mathrm{diff}}$ so the process will stop after finitely many steps, and we will obtain an element $\omega $ with norm less than $2\epsilon $ such that $\omega (p)\neq p$. Indeed, let $s = \displaystyle \mathop{\max }_{g_1, g_2\in \mathbb{S}_n, g_1\neq g_2}|Fix(g_1^{-1}g_2)|$. Then it suffices to take $k > s$.  $\square $
  
  \vspace{1cm} 
  
  \section {Weak Transitivity}

 In this section we prove Proposition \ref{thm:weakt} and then Proposition \ref{thm:second}. 
 
 \medskip
 
 First, we need the following lemma which follows immediately from the definition of $\Gamma _f$.
 
 \medskip
 
 \begin{lem} \label{lem:cute} Let $\Gamma $ be a finitely generated irreducible subgroup of the class $\mathbf{\Phi }^{\mathrm{diff}}$, $f$ be the biggest generator of $\Gamma $ with at least one fixed point, $z = \min Fix(f)$.  Then there exists $\epsilon > 0$ such that if $0 < a < b < \epsilon , f^p(a) < b$ for some $p\geq 4$, then, for all $\omega \in \Gamma _f$ with $\omega (z) > z$, if $n$  is sufficiently big, then $f^{p-4}(\omega f^{-n}\omega ^{-1}(a)) < \omega f^{-n}\omega ^{-1}(b)$.     
 \end{lem}
 
 \medskip 
 
 {\bf Proof.} Indeed, for sufficiently small $\epsilon $, if $0 <a < b < \epsilon < z$, then we have  $$f^{p-4}(\omega f^{-n}\omega ^{-1}(a)) <  f^{p-4+2-n}(a) = f^{-n-2}(f^p(a)) <  f^{-n-2}(b) < \omega f^{-n}\omega ^{-1} (b).$$  \ $\square $
 
 \medskip
 
 In the application of the lemma, our $\omega $ will be given in advance so we can use the following simpler version of the lemma as well:    {\em Let $\Gamma $ be a finitely generated irreducible subgroup of the class $\mathbf{\Phi }^{\mathrm{diff}}$, $f$ be the biggest generator of $\Gamma $ and $\omega \in \Gamma _f$ . Then there exists $\epsilon > 0$ such that if $0 < a < b < \epsilon , f^p(a) < b$ for some $p\geq 4$, then $f^{p-4}(\omega f^{-n}\omega ^{-1}(a)) < \omega f^{-n}\omega ^{-1}(b)$ for all sufficiently big $n$. } 
 
 \medskip

 {\bf Proof of Proposition \ref{thm:weakt}.} Let $\Gamma $ be an irreducible non-affine subgroup of $\mathbf{\Phi }^{\mathrm{diff}}, g\in \Gamma $, and $0 < p_1 < \dots < p_k < 1$. We may assume that $\Gamma $ is finitely generated. Then, let $f$ be the biggest generator of $\Gamma $. Without loss of generality we may also assume that $f$ has at least one fixed point and $\min Fix(f) < p_1$, moreover,  $f(x) > x$ and $g(x) > x$  for all sufficiently small positive $x$.  
 
 \medskip
 
 By dynamical 1-transitivity (i.e. Proposition \ref{thm:first}), there exist elements $\omega _1, \dots , \omega _k\in \Gamma _f$ such that $\min Fix(\eta _i) \in (p_{i-1}, p_{i}), 1\leq i\leq k$ where $\eta _i = \omega _i\dots \omega _1f\omega _1^{-1}\dots \omega _i^{-1}$.  (We also let $p_0 = 0, p_{k+1} = 1$.)
 
 \medskip

 We can choose $q\geq 16k$ such that $f^{[q/2]}(x) > g(x)$ for all $x\in (0, \frac{1}{2}\min Fix(f))$. By Lemma \ref{lem:cute}, for sufficiently big $n\geq 1$ and sufficiently small $\epsilon >0$,  we have  $$f^{4q-4}(\eta _1^{-n}(x)) < \eta _1^{-n}(y),\ \mathrm{for \ all} \ x, y\in (0, \epsilon ) \ \mathrm{where} \ f^{4q}(x) < y.$$

 Then, applying Lemma \ref{lem:cute} inductively, for sufficiently small $\epsilon > 0$ and sufficiently big $n$, we obtain that
 
  $$f^{q-16}(\eta _2^{-n}(p_1)) < \eta _2^{-n}(p_2), $$ \ $$\ldots $$ \ $$f^{q-8i}(\eta _i^{-n}(p_j)) < \eta _i^{-n}(p_{j+1}), 1\leq j\leq i-1$$    \ $$\dots $$ \ $$f^{q-8k}(\eta _k^{-n}(p_j)) < \eta _k^{-n}(p_{j+1}), 1\leq j\leq k-1.$$ 
 
 \medskip 
 
 Indeed,  for $2\leq i\leq k$, if $j = i-1$, then  the inequality follows from the fact that   $\min Fix(\eta _i) \in (p_{i-1}, p_{i})$. For $j\leq i-2$, by inductive hypothesis, we know that for sufficiently big $n$, $$f^{q-8(i-1)}(\eta _{i-1}^{-n}(p_j)) < \eta _{i-1}^{-n}(p_{j+1}).$$ Then for sufficiently big $n$, we have $\eta _{i-1}^{q-8i+6}(\eta _{i-1}^{-n}(p_j)) < \eta _{i-1}^{-n}(p_{j+1})$. This yields  an inequality $\eta _{i-1}^{q-8i+6}(p_j) < p_{j+1}$. By Lemma \ref{lem:cute} (and also because $\omega _i(\min Fix  \eta _{i-1}) > \min Fix  \eta _{i-1})$,
    
     $$\eta _{i-1}^{q-8i+2}(\omega _i\eta _{i-1}^{-n}\omega _i^{-1}(p_j)) < \omega _i\eta _{i-1}^{-n}\omega _i^{-1}(p_{j+1}).$$

\medskip 

     Then for sufficiently big $n$, $$f^{q-8i}(\omega _i\eta _{i-1}^{-n}\omega _i^{-1}(p_j)) < \omega _i\eta _{i-1}^{-n}\omega _i^{-1}(p_{j+1}).$$ 

     \medskip 
     
     Then $f^{q-8i}(\eta _{i}^{-n}(p_j)) < \eta _{i}^{-n}(p_{j+1})$. This completes the step of the induction.

     \medskip 
     
     Thus we proved the inequality $f^{q-4k}(\eta _k^{-n}(p_j)) < \eta _k^{-n}(p_{j+1}), 1\leq j\leq k-1$. Now,  it suffices to take $\gamma = \eta _k^{-n}$, and we obtain that $$\gamma (p_{j-1}) < f^q(\gamma (p_j)) < \gamma (p_{j+1}), 1\leq j\leq k-1.$$  Then, we have $$\gamma (p_{j-1}) < g(\gamma (p_j)) < \gamma (p_{j+1}), 1\leq j\leq k-1.$$ 
    $\square $ 
 
 \bigskip
 
 Now we are ready to prove Proposition \ref{thm:second}. Let $\Gamma $ be a finitely generated irreducible non-affine subgroup of $\mathrm{Diff}_{+}(I)$. 
 
 \medskip
 
 Let $S = \{f_1, \dots , f_s\}$ be a generating set of $\Gamma $ such that $S\cap S^{-1} = \emptyset $ (in particular, $S$ does not contain the identity element), $f$ be the biggest generator of $S\cup S^{-1}$, and let also $m\geq 10s$. By the definition of the order, there exists $\epsilon > 0$ such that $f(x) \geq \xi (x)$ for all $x\in (0,\epsilon ), \xi \in S\cup S^{-1}$.  
 
 \medskip
 
 By Proposition \ref{thm:third} we can choose $g\in \Gamma $ such that $|Fix(g)| > 8m$. Since $\Gamma $ is irreducible, by weak transitivity, we can find $h\in \Gamma $ and $\epsilon > 0$ such that the following conditions hold:
 
 \medskip

  (i) $\gamma (\mathrm{Fix}(hgh^{-1})) \subset (0, \epsilon )$ for all $\gamma \in B_{2m}(1)$,
 
 \medskip
 
  (ii) if $Fix(hgh^{-1}) = \{p_1, \dots , p_k\}$, then $p_{i-1} < f^{2m} (p_i) < p_{i+1}, 2\leq i\leq k-1$.
   
 \medskip
 
  (iii) $Fix(f_i)\cap Fix(hgh^{-1}) = \emptyset $, for all $i\in \{1, \dots , s\}$.   
  
  \medskip
 
  Now, let $\theta = hgh^{-1}$ and $$S_n = \{\theta , \theta ^{mn}f_1\theta ^{mn},  \theta ^{2mn}f_2\theta ^{2mn}, \dots ,  \theta ^{smn}f_s\theta ^{smn}\}.$$ Let also $$\delta = \frac{1}{5}\min \{p_{i+1}-p_i \ | \ 1\leq i\leq k-1\} \ \mathrm{and} \  V = \displaystyle \mathop{\sqcup} _{2\leq i\leq k-1}(p_i-\delta , p_i + \delta ).$$ 
  
  \medskip
  
  We now let $x = \frac{p_j+p_{j+1}}{2}$ where $j = [\frac{k}{2}]$. Then $x\notin V$, and for a sufficiently big $n$, by a standard ping-pong argument, for any non-trivial word $W$ in the alphahbet $S_n$ of length at most $m$, we obtain that $W(x)\in V$, hence $W(x) \neq x$.  
 
 \medskip
 
 Since $m$ is arbitrary, we conclude that $girth(\Gamma ) = \infty $.


\vspace{1cm}
  
  
  
  \begin{thebibliography}{99}
  
  \bibitem{A1} Akhmedov A. \ A weak Zassenhaus lemma for subgroups of Diff(I). {\em Algebraic and Geometric Topology}. {\bf vol.14} (2014) 539-550.
  
  \bibitem{A2} Akhmedov A. \ Extension of H\"older's Theorem in $\mathrm{Diff}_{+}^{1+\epsilon }(I)$.  \ {\em Ergodic Theory and Dynamical Systems}, {\bf  vol.36}, no 5, (2016), 1343-1353. 
  
  \bibitem{A3} Akhmedov, A. \ On groups of homeomorphisms of the interval with finitely many fixed points. \  Preprint. https://arxiv.org/abs/1503.03850
  
  \bibitem{A4} Akhmedov, A. \ Girth Alternative for subgroups of $PL_{o}(I)$. \ Preprint. http://arxiv.org/pdf/1105.4908.pdf
  
  \bibitem{A5} Akhmedov, A. \ On the height of subgroups of $\mathrm{Homeo}_{+}(I)$, {\em Journal of Group Theory}, {\bf 18}, no.1, (2015) 93-108.
  
   \bibitem{A6} Akhmedov, A.\  On the girth of finitely generated groups, Journal of Algebra 268 (2003), 198-208.
   
    \bibitem{A7} Akhmedov, A.\ The girth of groups satisfying Tits Alternative, Journal of Algebra 287 (2005), 275–282.
    
     \bibitem{A8} Akhmedov, A.\ Questions and remarks on discrete and dense subgroups of Diff(I), Journal of Topology and Analysis, vol. 6, no. 4, (2014), 557-571.
     
       \bibitem{A9} Akhmedov, A.\ On free discrete subgroups of Diff(I).  Algebraic and Geometric Topology, vol.10, no.4, (2010) 2409-2418. 
  
   \bibitem{FF} B.Farb, J.Franks, Groups of homeomorphisms of one-manifolds II: Extension of H\"older's Theorem. Trans. Amer. Math. Soc. 355 (2003) no.11, 4385-4396.
   
   \bibitem{MN} Bonatti, C, Monteverde, I., Navas, A., Rivas, C. \ Rigidity for $C^1$ actions on the interval arising from hyperbolicity I: solvable groups. Mathematische Zeitschrift,  286 (2017), 919–949.
   
     \bibitem{M} M. Muller. \ Sur l’approximation et l’instabilit´e des feuilletages. Unpublished text (1982).
  
  \bibitem{N1} Navas, A. \ Groups, Orders and Laws. \ {\em Groups, Geometry, Dynamics.}   8 (2014), 863–882, \  http://arxiv.org/abs/1405.0912 
  
  \bibitem{N2} Navas, A. \ Groups of Circle Diffemorphisms. \ Chicago Lectures in Mathematics, 2011. {\em http://arxiv.org/pdf/0607481} 
  
   \bibitem{S} Schleimer, S. \ On the girth of groups. Preprint.\
  
  \bibitem{T} T. Tsuboi.\  $G_1$-structures avec une seule feuille. Asterisque. 116 (1984), 222-234.
  
 
  
  
  \end{thebibliography}
 \end{document}